
\magnification=1200
\hoffset3cm

\input amstex

\catcode`\@=11
\font\myfont=cmss10
\def\myifdefined#1#2#3{%
  \expandafter\ifx\csname #1\endcsname \relax#2\else#3
    \fi}
\newcount\fornumber\newcount\artnumber\newcount\tnumber
\newcount\secnumber
\newcount\constnumber
\newcount\theonumber
\def\Ref#1{%
  \expandafter\ifx\csname mcw#1\endcsname \relax
    \warning{\string\Ref\string{\string#1\string}?}%
    \hbox{$???$}%
  \else \csname mcw#1\endcsname \fi}
\def\Refpage#1{%
  \expandafter\ifx\csname dw#1\endcsname \relax
    \warning{\string\Refpage\string{\string#1\string}?}%
    \hbox{$???$}%
  \else \csname dw#1\endcsname \fi}

\def\warning#1{\immediate\write16{%
            -- warning -- #1}}
\def\CrossWord#1#2#3{%
  \def\x{}%
  \def\y{#2}%
  \ifx \x\y \def\z{#3}\else
            \def\z{#2}\fi
  \expandafter\edef\csname mcw#1\endcsname{\z}
\expandafter\edef\csname dw#1\endcsname{#3}}
\def\Talg#1#2{\begingroup
  \edef\mmhead{\string\CrossWord{#1}{#2}}%
  \def\writeref{\write\refout}%
  \expandafter \expandafter \expandafter
  \writeref\expandafter{\mmhead{\the\pageno}}%
\endgroup}

\openin15=\jobname.ref
\ifeof15 \immediate\write16{No file \jobname.ref}%
\else      \input \jobname.ref \fi \closein15
\newwrite\refout
\openout\refout=\jobname.ref

\def\dfs#1{\myifdefined{mcwx#1}{}{\warning{multiply defined label #1}}
     \global\advance\secnumber by 1\global\fornumber=0\global\tnumber=0
      \the\secnumber.\if\the\theonumber0
      {\myfont #1}%
\else\relax\fi
     \Talg{#1}{\the\secnumber}
 \expandafter\gdef\csname mcwx#1\endcsname{a}\ignorespaces
}%

\def\dff#1{\myifdefined{mcwx#1}{}{\warning{multiply defined label #1}}
     \global\advance\fornumber by 1 \if\the\theonumber0 \text{\myfont #1}%
\else\relax\fi
\unskip (\the\secnumber.\the\fornumber)
\enspace
     \Talg{#1}{(\the\secnumber.\the\fornumber)}
  \expandafter\gdef\csname mcwx#1\endcsname{a} \ignorespaces
}%
\def\dfc#1{\myifdefined{mcwx#1}{}{\warning{multiply defined label #1}}
     \global\advance\constnumber by 1 \if\the\theonumber0 \text{\myfont #1}%
\else\relax\fi
\unskip \the\constnumber
     \Talg{#1}{\the\constnumber}
  \expandafter\gdef\csname mcwx#1\endcsname{a} \ignorespaces
}%


\def\dft#1{\myifdefined{mcwx#1}{}{\warning{multiply defined label #1}}%
     \global\advance\tnumber by 1%
       \if\the\theonumber0%
 \ignorespaces    {\myfont #1 }%
\else\relax\fi%
 \the\secnumber.\the\tnumber
    \Talg{#1}{\the\secnumber.\the\tnumber}
  \expandafter\gdef\csname mcwx#1\endcsname{a} \unskip\ignorespaces
}%
\def\dfro#1{\myifdefined{mcwx#1}{}{\warning{multiply defined label #1}}%
     \if\the\theonumber0%
      {\par\noindent\llap{\myfont #1}\par}%
\else\relax\fi%
     \Talg{#1}{(\the\rostercount@)}%
  \expandafter\gdef\csname mcwx#1\endcsname{a} \unskip\ignorespaces
}%
\def\dfa#1{\myifdefined{mcwx#1}{}{\warning{multiply defined label #1}}
     \global\advance\artnumber by 1
      \unskip      \the\artnumber%
     \Talg{art.#1}{\the\artnumber}
     \expandafter\gdef\csname mcwx#1\endcsname{a} \unskip\ignorespaces
 }%

\def\rf#1{\ifmmode\Ref{#1}\else $\Ref{#1}$\fi}

\def\donotshowtheoremlabels{\theonumber=1}

\def\rfa#1{\Ref{art.#1}}

\csname amsppt.sty\endcsname

\documentstyle{amsppt}

\vsize 18.8 cm


\def\df#1{\leqno{\myifdefined{mcwx#1}{}{\warning{multiply
defined label #1}}\unskip
     \global\advance\fornumber by 1 \if\the\formunumber0 {\rm ref\{#1\}}%
\else\relax\fi
     \unskip (\the\secnumber.\the\fornumber)
    \Talg{#1}{(\the\secnumber.\the\fornumber)}\unskip
   \expandafter\gdef\csname mcwx#1\endcsname{a}\ignorespaces
}}%

\def\dfr#1#2{\myifdefined{mcwx#1}{}{\warning{multiply
defined label #1}}
     \Talg{#1}{#2}
  \expandafter\gdef\csname mcwx#1\endcsname{a} \unskip\ignorespaces
}%

\newdimen\mydim
\def\mycenter#1{\mydim=\hsize \advance \mydim by -50pt
\vtop{\hsize=\mydim
 \noindent\ignorespaces  #1 }}

\def\proof{\demo{Proof}} \def\eproof{\qed\enddemo}
\def\co{\colon}

\def\N{{\Bbb N}}
\def\Z{{\Bbb Z}}

\def\R{{\Bbb R}}

\def\D{\roman{d}}
\def\D{\roman{d}}
\def\eps{\varepsilon}
\def\is#1..#2..{\langle #1,#2\rangle}
\def\cici#1..#2..{\mathinner{[{#1}\ldotp\ldotp{#2}]}}
\def\ci#1..{\mathinner{[{#1}]}}
\def\oi#1..{\mathinner{]{#1}[}}
\def\ro#1..{\mathinner{[{#1}[}}
\def\lo#1..{\mathinner{]{#1}]}}
\def\ici#1..{\mathinner{[\![{#1}]\!]}} \def\Om{\Omega}
\donotshowtheoremlabels 
\nologo
\hoffset1.5cm

\topmatter
\def\mytime{\the\day.\the\month.\the\year-\the\time}
\title  Attractors for singularly
perturbed hyperbolic equations on unbounded domains
\endtitle

\author  Martino Prizzi --- Krzysztof P.
Rybakowski
\endauthor
\leftheadtext{ M. Prizzi --- K. P.
Rybakowski }
\rightheadtext{Damped hyperbolic equations  }

\address Martino Prizzi, Universit\`a degli Studi di
Trieste, Dipartimento di Matematica e Informatica, Via
Valerio, 12, 34127 Trieste, ITALY
\endaddress
\email
prizzi\@ dsm.univ.trieste.it
\endemail

\address Krzysztof P. Rybakowski, Universit\"at Rostock,
Institut f\"ur Mathematik, Universit\"atsplatz 1, 18055
Rostock, GERMANY
\endaddress
\email
krzysztof.rybakowski\@ mathematik.uni-rostock.de
\endemail

\abstract
For an arbitrary unbounded domain $\Omega\subset\R^3$ and for $\eps>0$, we consider the
damped hyperbolic equations
$$\aligned \eps
u_{tt}+ u_t+\beta(x)u- \sum_{ij}(a_{ij}(x)
u_{x_j})_{x_i}&=f(x,u),\quad x\in \Omega,\,t\in\ro0,\infty..,\\
u(x,t)&=0,\quad x\in \partial \Omega,\, t\in\ro0,\infty...
\endaligned\leqno{(H_\eps)}$$
and their singular limit as $\eps\to0$, i.e. the parabolic equation
$$\aligned
u_t+\beta(x)u- \sum_{ij}(a_{ij}(x)
u_{x_j})_{x_i}&=f(x,u),\quad x\in \Omega,\,t\in\ro0,\infty..,\\
u(x,t)&=0,\quad x\in \partial \Omega,\, t\in\ro0,\infty...
\endaligned\leqno{(P)}$$
Under suitable assumptions, $(H_\eps)$ possesses a compact global attractor $\Cal A_\eps$
in the phase space $H^1_0(\Omega)\times L^2(\Omega)$, while $(P)$ possesses a compact global attractor
$\widetilde{\Cal A_0}$ in the phase space $H^1_0(\Omega)$, which can be embedded into a compact set
${\Cal A_0}\subset H^1_0(\Omega)\times L^2(\Omega)$. We show that, as $\eps\to0$, the family
$({\Cal A_\eps})_{\eps\in[0,\infty[}$ is upper semicontinuous with respect to the topology of
$H^1_0(\Omega)\times H^{-1}(\Omega)$. We thus extend a well known result by Hale and Raugel in three
directions: first, we allow $f$ to have critical growth; second, we let $\Omega$ be unbounded; last,
we do not make any smoothness assumption on $\partial\Omega$, $\beta(\cdot)$, $a_{ij}(\cdot)$ and $f(\cdot,u)$.
\endabstract

\endtopmatter

\document \hoffset0cm

\head
\dfs{sec:intr}Introduction
\endhead

In their paper \cite{\rfa{HR}} Hale and Raugel considered the
damped hyperbolic equations
$$\aligned \eps
u_{tt}+ u_t-\Delta u&
=f(u)+g(x),\quad x\in \Omega,\,t\in\ro0,\infty..,\\
u(x,t)&=0,\quad x\in \partial \Omega,\, t\in\ro0,\infty...
\endaligned$$
and their singular limit as $\eps\to0$, i.e. the parabolic equation
$$\aligned
u_t-\Delta u&=f(u)+g(x),\quad x\in \Omega,\,t\in\ro0,\infty..,\\
u(x,t)&=0,\quad x\in \partial \Omega,\, t\in\ro0,\infty...
\endaligned$$
In \cite{\rfa{HR}} the set $\Omega$ is a bounded smooth domain or a convex polyhedron, $\eps$ is a positive constant,
$g\in L^2(\Omega)$ and $f$ is a $C^2$ function of subcritical growth such that
$$\limsup_{|u|\to\infty}{f(u)\over u}\leq 0.$$
Under these assumptions, for any fixed $\eps>0$ the corresponding hyperbolic equation generates a global
semiflow which possesses a compact global attractor $\Cal A_\eps$
in the phase space $H^1_0(\Omega)\times L^2(\Omega)$ (see \cite{\rfa{BV},\rfa{GT},\rfa{Ha}}). Moreover, the
limiting parabolic equation generates a  global semiflow which possesses a compact global attractor
$\widetilde{\Cal A_0}$ in the phase space $H^1_0(\Omega)$ (see \cite{\rfa{DCh},\rfa{Ha}}). Due to the smoothing
effect of parabolic equations, it turns out that $\widetilde{\Cal A_0}$ is actually a compact subset of
$H^2(\Omega)$. Hence one can define the set
$${\Cal A_0}=\{(u,\Delta u+f(u)+g)\mid u\in{\Cal A_0} \},$$
which is a compact subset of $H^1_0(\Omega)\times L^2(\Omega)$. Hale and Raugel proved that
the family
$({\Cal A_\eps})_{\eps\in[0,\infty[}$ is upper semicontinuous with respect to the topology of
$H^1_0(\Omega)\times L^{2}(\Omega)$, i.e.
$$\lim_{\eps\to 0^+}\sup_{y\in\Cal A_\eps}\inf_{z\in \Cal A_0}|y-z|_{H^1_0\times L^{2}}=0.$$

In this paper we extend the result of Hale and Raugel in three
directions: firstly, we allow $f$ to have critical growth; secondly, we let $\Omega$ be unbounded; thirdly,
we replace $f(u)+g(x)$ by $f(x,u)$ and $-\Delta$ by $\beta(x)u- \sum_{ij}(a_{ij}(x)
u_{x_j})_{x_i}$, without any smoothness assumption on $\partial\Omega$,
$\beta(\cdot)$, $a_{ij}(\cdot)$ and $f(\cdot,u)$.

In \cite{\rfa{HR}} the proof of the main result relies on some uniform $(H^2\times H^1)$-estimates
for the attractors $\Cal A_\eps$, combined with the compactness of the Sobolev embedding
$H^1_0(\Omega)\subset L^2(\Omega)$. The uniform $(H^2\times H^1)$-estimates are obtained through a
bootstrapping argument originally due to Haraux \cite{\rfa{Har}}. Such argument works only if $f$ is
subcritical, and if $\Omega$ is such that the domain of the $L^2(\Omega)$-realization of $-\Delta$ is
$H^2(\Omega)\cap H^1_0(\Omega)$ (e.g. if $\Omega$ is a convex polyhedron).

A different bootstrapping argument was proposed by Grasselli and
Pata in \cite{\rfa{GraPa1},\rfa{GraPa2}}. Their argument also
works  in the critical case, and is based on certain a-priori
estimates that can be obtained ``within an appropriate Galerkin
approximation scheme".  Here, ``appropriate" means ``on a basis of
eigenfunctions of $-\Delta$".  Therefore, their approach cannot be
used in the case of an unbounded domain $\Omega$. More recently,
in \cite{\rfa{PaZe}} Pata and Zelik obtained $(H^2\times
H^1)$-estimates for ${\Cal A}_\eps$ without using bootstrapping
arguments, but again their a-priori estimates are obtained
``within an appropriate Galerkin approximation scheme". We point
out that also in \cite{\rfa{GraPa1},\rfa{GraPa2},\rfa{PaZe}}
$\Omega$ must  have the property that the domain of the
$L^2(\Omega)$-realization of $-\Delta$ is $H^2(\Omega)\cap
H^1_0(\Omega)$. Moreover, the Nemitski operator associated with
$f$ must be Lipschitz continuous from $H^2(\Omega)\cap
H^1_0(\Omega)$ to $H^1(\Omega)$ in \cite{\rfa{PaZe}} and from
$D((-\Delta)^{(\alpha+1)/2})$ to $D((-\Delta)^{\alpha/2})$ for all
$0\leq\alpha\leq 1$ in \cite{\rfa{GraPa1},\rfa{GraPa2}}.
Therefore, if one wants to replace $f(u)+g(x)$ by $f(x,u)$, one
needs to impose severe smoothness conditions on $f(x,u)$ with
respect to the space variable $x$.

If $\Omega$ is unbounded, the embedding $H^1_0(\Omega)\subset L^2(\Omega)$ is no longer compact, and this poses
some additional difficulties even  for the   existence proof of the attractors ${\Cal A_\eps}$.
In \cite{\rfa{F},\rfa{F1}}, Feireisl circumvented these difficulties by decomposing any solution $u(t,x)$
into the sum $u_1(t,x)+u_2(t,x)$ of two functions, such that $u_1(t,\cdot)$ is asymptotically small,
and $u_2(t,\cdot)$ has a compact support which propagates with speed $1/\eps^2$. As $\eps\to 0$, the speed of
propagation tends to infinity, and, indeed, the estimates obtained by Feireisl are not uniform with respect to
$\eps$. It is therefore apparent that, if one wants to pass to the limit as $\eps\to0$, a different approach is
needed.

In our previous paper \cite{\rfa{PR}} we proved the existence of compact global attractors for damped
hyperbolic equations in unbounded domains using the method of tail-estimates
(introduced by Wang in \cite{\rfa{W}} for parabolic equations), combined with an argument due to Ball
\cite{\rfa{Ba}} and elaborated by Raugel in \cite{\rfa{Ra}}. Here we exploit the same techniques to establish
an upper semicontinuity result similar to that of Hale and Raugel, when $\Omega$ is an unbounded domain and $f$
is critical. Our arguments do not rely on $(H^2\times H^1)$-estimates for the attractors $\Cal A_\eps$.
Therefore they also apply  to the case of an open set $\Omega$ for which
the domain of the $L^2(\Omega)$-realization of $-\Delta$ is not
$H^2(\Omega)\cap H^1_0(\Omega)$ (e.g. if $\Omega$ is the exterior of a convex polyhedron).

Before we describe in detail our assumptions and our results, we
need to introduce some notation. In this paper, $N=3$ and $\Omega$
is an arbitrary open subset of $\R^N$, bounded or not. For $a$ and
$b\in\Z$ we write $\cici a..b..$ to denote the set of all $m\in\Z$
with $a\le m\le b$. Given a subset $S$ of $\R^N$ and a function
$v\co S\to \R$ we denote by  $\tilde v\co\R^N\to \R$  the trivial
extension of $v$ defined by $\tilde v(x)=0$ for $x\in
\R^N\setminus S$. Given  a function $g\co \Omega\times \R\to \R$,
we denote by $\hat g$  the Nemitski operator which
associates with every function $u\co \Omega\to \R$ the function
$\hat g(u)\co \Omega\to \R$ defined by
 $$
  \hat g(u)(x)= g(x,u(x)),\quad x\in \Omega.
 $$
Unless specified otherwise, given $k\in \N$ and functions $g$, $h\co\Omega\to \R^k$ we write $$\langle g,h\rangle:=\int_\Omega \sum_{m=1}^k g_m(x)h_m(x)\,\D x,$$ whenever the integral on the right-hand side makes sense.

If $I\subset \R$, $Y$ and $X$
are normed spaces with $Y\subset X$ and if $u\co I\to Y$ is a
function which is differentiable as a function into $X$ then we
denote its $X$-valued derivative by $\partial(u;X)$. Similarly, if
$X$ is a Banach space and $u\co I\to X$ is integrable as a
function into $X$, then we denote its $X$-valued integral by $\int
_I(u(t);X)\,\D t$.

\proclaim{Assumption~\dft{290506-1840}}\roster
\item\dfro{290506-1841} $a_0$,
$a_1\in\oi0,\infty..$ are constants and $a_{ij}\co
\Omega\to \R$, $i$, $j\in\cici1..N..$ are functions in
$L^\infty(\Omega)$ such that $a_{ij}=a_{ji}$, $i$,
$j\in\cici1..N..$, and for every $\xi\in\R^N$ and a.e.
$x\in\Omega$, $a_0|\xi|^2\le \sum_{i,j=1}^N
a_{ij}(x)\xi_i\xi_j\le a_1|\xi|^2 $.
$A(x):=(a_{ij}(x))_{i,j=1}^N$,
$x\in\Omega$.\item\dfro{290506-1842} $\beta\co
\Omega\to \R$ is a measurable function with the property
that \itemitem{$(i)$}\dfr{290506-1843}{(2i)} for every
$\overline \eps\in \oi0,\infty..$ there is a
$C_{\overline\eps}\in \ro0,\infty..$ with
$\bigl||\beta|^{1/2} u\bigr|_{L^2}^2\le \overline \eps
|u|_{H^1}^2+C_{\overline\eps}|u|_{L^2}^2 $ for all $u\in
H^1_0(\Omega)$;\itemitem{$(ii)$} $\lambda_1:=\inf\{\,\is
A\nabla u..\nabla u.. +\langle \beta  u,u\rangle \mid u\in
H^1_0(\Om),\,|u|_{L^2}=1\,\}>0$.\endroster\endproclaim

\proclaim{Assumption~\dft{140506-1137}}
\roster\item\dfro{290506-1844}
$f\co \Omega\times\R\to \R$ is such that, for every $u\in\R$, $f(\cdot, u)$ is (Lebesgue-)measurable, $f(\cdot,0)\in L^2(\Omega)$ and for a.e. $x\in\Omega$, $f(x,\cdot)$ is of class $C^2$ and such that $\partial_u f(\cdot,0)\in L^\infty(\Omega)$ and $|\partial_{uu}f(x,u)|\le \overline C(1+|u|)$ for some constant $\overline C\in\ro0,\infty..$, every $u\in \R$ and a.e. $x\in\Omega$;
 \item\dfro{290506-1845} $f(x,u)u-\overline\mu F(x,u)\le
c(x)$ and $F(x,u)\le c(x)$ for a.e.  $x\in \Om$ and every
$u\in\R$. Here, $c\in L^2(\Omega)$ is a given function, $\overline\mu\in\oi0,\infty..$ is a constant and $F\co \Om\times \R\to \R$ is defined, for $(x,u)\in\R$, by $$F(x,u)=\int_0^u
f(x,s)\,\D s,
 $$ whenever $f(x,\cdot)\co \R\to \R$ is continuous and $F(x,u)=0$
 otherwise.
\endroster
\endproclaim
Note that Assumptions~\rf{290506-1840} and~\rf{140506-1137} imply the hypotheses of~\cite{\rfa{PR}}.

Let $D(\bold B_{\eps})$ be the set of all $(u,v)\in
H^1_0(\Omega)\times L^2(\Omega)$ such that $v\in H^1_0(\Om)$ and
$-\beta u+\sum_{ij}(a_{ij} u_{x_j})_{x_i}$ (in the distributional
sense) lies in $L^2(\Om)$. It turns out that the operator $$\bold
B_\eps(u,v)=(-v,(1/\eps) v+(1/\eps)\beta
u-(1/\eps)\sum_{ij}(a_{ij} u_{x_j})_{x_i}),\quad (u,v)\in D(\bold
B_\eps)$$  is the generator of a $(C_0)$-semigroup $e^{-\bold
B_\eps t}$, $t\in\ro0,\infty..$ on $H^1_0(\Omega)\times
L^2(\Omega)$. Moreover, the Nemitski operator $\hat f$ is a
Lipschitzian map of $H^1_0(\Omega)$ to $L^2(\Omega)$. Results in
\cite{\rfa{CH}} then imply that  the hyperbolic boundary value
problem $$\aligned \eps u_{tt}+ u_t+\beta(x)u- \sum_{ij}(a_{ij}(x)
u_{x_j})_{x_i}&=f(x,u),\quad x\in \Omega,\,t\in\ro0,\infty..,\\
u(x,t)&=0,\quad x\in \partial \Omega,\, t\in\ro0,\infty..
\endaligned$$
with Cauchy data at $t=0$ has a unique (mild) solution $z(t)=(u(t),v(t))$ in $H^1_0(\Omega)\times L^2(\Omega)$, given by the ``variation-of-constants" formula
$$
z(t)=e^{-\bold B_\eps t}z(0)+\int_0^te^{-\bold B_\eps(t-s)}(0,(1/\eps)\hat f(u(s)))\,\D s.
$$
For $\eps\in\oi0,\infty..$ we define $\pi_\eps$ to be the local semiflow on $H^1_0(\Omega)\times L^2(\Omega)$ generated by the (mild) solutions of this hyperbolic boundary value problem.
We can summarize the results of \cite{\rfa{PR}} in the following:

\proclaim{Theorem~\dft{230107-1611}} Under
Assumptions~\rf{290506-1840} and~\rf{140506-1137},
$\pi_\eps$ is a global semiflow and it has a global attractor $\Cal A_\eps$.
\endproclaim

Analogously, consider the parabolic boundary value problem
$$\aligned
u_t+\beta(x)u- \sum_{ij}(a_{ij}(x)
u_{x_j})_{x_i}&=f(x,u),\quad x\in \Omega,\,t\in\ro0,\infty..,\\
u(x,t)&=0,\quad x\in \partial \Omega,\, t\in\ro0,\infty..
\endaligned$$
with Cauchy data at $t=0$.
Letting $\bold A$ denote the sectorial
operator on
$L^2(\Omega)$ defined by the differential operator
$u\mapsto \beta u- \sum_{ij}(a_{ij}u_{x_j})_{x_i}$, we
have that
$D(\bold A)$ is the set of all $u\in H^1_0(\Omega)$ such that the distribution
$\beta u-\sum_{ij}(a_{ij}u_{x_j})_{x_i}$ lies in $L^2(\Omega)$. Again, the Cauchy problem has a unique (mild) solution $u(t)$ in $H^1_0(\Omega)$, given by the ``variation-of-constants" formula
$$
u(t)=e^{-\bold A t}u(0)+\int_0^te^{-\bold A(t-s)}\hat f(u(s))\,\D s.
$$
Let $\widetilde\pi$ be the local semiflow on $H^1_0(\Omega)$ generated by the (mild) solutions of this parabolic boundary value problem.
Results in~\cite{\rfa{PR3}} imply that $\widetilde \pi$ is a global semiflow and has a global
attractor~$\widetilde{\Cal A}$ (see also \cite{\rfa{ACDR}}). Moreover, it is proved in~\cite{\rfa{PR3}} that $\widetilde{\Cal A}\subset D({\bold A})$ and $\widetilde{\Cal A}$ is compact in $D({\bold A})$ endowed with the graph norm.

 Let $\Gamma\co D({\bold A})\to H^1_0(\Omega)\times L^2(\Omega)$ be defined by $\Gamma(u)=(u,{\bold A}u+\widehat f(u))$.  Set $\Cal A_0:=\Gamma(\widetilde{\Cal A})$. Then we have the following main result of this paper:
\proclaim{Theorem~\dft{300506-0938}} The family $(\Cal A_\eps)_{\eps\in\ro0,\infty..}$ is upper semicontinuous at $\eps=0$ with respect to the topology of $H^1_0(\Omega)\times H^{-1}(\Omega)$, i.e.
$$\lim_{\eps\to 0^+}\sup_{y\in\Cal A_\eps}\inf_{z\in \Cal A_0}|y-z|_{H^1_0\times H^{-1}}=0.$$
\endproclaim
Actually a stronger result is established in Theorem~\rf{151106-1022} below.

\head
\dfs{300506-1021} Preliminaries
\endhead

In this section we collect a few preliminary results.
We begin with an abstract lemma established in~\cite{\rfa{PR3}}:
\proclaim{Lemma~\dft{220506-0819}}
Suppose $(Y,\langle\cdot,\cdot\rangle_{Y})$ and $(X,\langle\cdot,\cdot\rangle_{X})$ are (real or complex)
Hilbert spaces such that $Y\subset X$, $Y$ is dense in $(X,\langle\cdot,\cdot\rangle_{X})$ and the inclusion
$(Y,\langle\cdot,\cdot\rangle_{Y})\to(X,\langle\cdot,\cdot\rangle_{X})$
is continuous.
Then for every $u\in X$ there exists a unique $w_u\in Y$ such that
$$\langle v,w_u\rangle_Y=\langle v,u\rangle_X\text{ for all $v\in Y$.}$$
The map $B\co X\to X$, $u\mapsto w_u$ is linear, symmetric and positive. Let $B^{1/2}$ be a square root of $B$, i.e. $B^{1/2}\co X\to X$ linear, symmetric and $B^{1/2}\circ B^{1/2}=B$.
Then   $B$ and $B^{1/2}$ are injective and $R(B)$ is dense in $Y$. Set $X^{1/2}=X^{1/2}_B=R(B^{1/2})$ and $B^{-1/2}\co X^{1/2}\to X$ be the inverse of $B^{1/2}$. On $X^{1/2}$ the assignment $\langle u,v\rangle_{1/2}:=\langle B^{-1/2}u, B^{-1/2}v\rangle_X$ is a complete scalar product. We have $Y=X^{1/2}$ and $\langle \cdot,\cdot\rangle_{Y}=\langle \cdot,\cdot\rangle_{1/2}$.
\endproclaim
Now let $\bold A$ be the sectorial
operator on
$L^2(\Omega)$ defined by the differential operator
$u\mapsto \beta u- \sum_{ij}(a_{ij}u_{x_j})_{x_i}$.
Then ${\bold A}$ generates a family $X^\alpha=X^\alpha_{\bold A}$, $\alpha\in \R$, of fractional power spaces with $X^{-\alpha}$ being the dual of $X^\alpha$ for $\alpha\in\oi0,\infty..$. We write $$H_\alpha=X^{\alpha/2},\quad \alpha\in\R.$$ For $\alpha\in\R$ the operator $\bold A$ induces an operator ${\bold A}_\alpha\co H_\alpha\to H_{\alpha-2}$. In particular, $H_0=L^2(\Omega)$ and ${\bold A}={\bold A}_2$.

Note that, thanks to Assumption \rf{290506-1840}, the scalar product
$$\langle u,v\rangle_{H^1_0}=\langle A\nabla u,\nabla v\rangle+\langle \beta u,v\rangle, \quad u,v\in H^1_0(\Omega)$$
on $H^1_0(\Omega)$ is equivalent to the usual scalar product on $H^1_0(\Omega)$. Moreover,
$$\langle u,v\rangle_{H^1_0}=\langle {\bold A}_2 u, v\rangle, \quad u\in D(A_2),v\in H^1_0(\Omega).$$
\proclaim{Corollary~\dft{300506-1029}}
$H_1=H^1_0(\Omega)$ with equivalent norms. Consequently $H_{-1}=H^{-1}(\Omega)$ with equivalent norms.
\endproclaim

\proof Set $(X,\langle\cdot,\cdot\rangle_X)=(L^2(\Omega),\langle \cdot,\cdot\rangle)$ and $(Y,\langle\cdot,\cdot\rangle_Y)=(H^1_0(\Omega),\langle \cdot,\cdot\rangle_{H^1_0})$.
 Then $Y$ is dense in $X$ and the inclusion $Y\to X$ is continuous. Let $B_2\co X\to X$ be the inverse of ${\bold A}_2$. 
Then for all $u\in X$, $B_2u\in Y$
and for all $v\in Y$
$$\langle v,u \rangle_X=\langle v,B_2u\rangle_Y.$$ Thus $B_2=B$ where $B$ is as in Lemma~\rf{220506-0819}. Now the lemma implies the corollary.
\eproof
\proclaim{Corollary~\dft{230506-0733}}
The linear operator ${\bold A}_1\co H_1\to X:=H_{-1}$ is self-adjoint hence sectorial on $X$. Let $X_1^\alpha$, $\alpha\in\ro0,\infty..$, be the family of fractional powers generated by ${\bold A}_1$. Then $X^{1/2}=L^2(\Omega)$ with equivalent norms.
\endproclaim

\proof Set $(X,\langle\cdot,\cdot\rangle_X)=(H_{-1},\langle \cdot,\cdot\rangle_{H_{-1}})$ and $(Y,\langle\cdot,\cdot\rangle_Y)=(H_{0},\langle \cdot,\cdot\rangle_{H_{0}})$. Then $Y$ is dense in $X$ and the inclusion $Y\to X$ is continuous. Let $B_1\co X\to X$ be the inverse of ${\bold A}_1$.
Then for all $u\in X$, $B_1u\in Y$
and for all $v\in Y$
$$\langle v,u \rangle_X=\langle B_1v,B_1u\rangle_{H_1}=\langle v,B_1u\rangle_Y.$$ Thus $B_1=B$ where $B$ is as in Lemma~\rf{220506-0819}. Now the lemma implies the corollary.
\eproof

We end this section by quoting a result proved in \cite{\rfa{PR}}, which can be used to
rigorously justify formal differentiation of various
functionals along (mild) solutions of semilinear evolution
equations.
\proclaim{Theorem~\dft{230107-1639}}
 Let $Z$ be a Banach space and $B\co D(B)\subset Z\to Z$ the infinitesimal
 generator of a $(C_0)$-semigroup of linear operators $e^{-Bt}$ on $Z$, $t\in\ro0,\infty..$.
 Let $U$ be open in $Z$, $Y$ be a normed space and $V\co U\to Y$ be a function which, as a map from $Z$ to $Y$, is continuous at
 each point of $U$ and Fr\'echet differentiable at each
 point of $U\cap D(B)$. Moreover, let $W\co U\times Z\to Y$
 be a function which, as a map from $Z\times Z$ to $Y$, is
 continuous  and  such that $DV(z)(Bz+w)=W(z,w)$ for $z\in
 U\cap D(B)$ and $w\in Z$. Let $\tau\in \oi0,\infty..$ and
 $I:=\ci0,\tau..$. Let $\bar z\in U$, $g\co I\to Z$ be
 continuous and $z$ be a map from $I$ to $U$ such that
 $$
  z(t)=e^{-Bt}\bar z+\int_0^t e^{-B(t-s)}g(s)\,\D s, \quad t\in I.
 $$
 Then the map $V\circ z\co I\to Y$ is differentiable
 and
 $$
  (V\circ z)'(t)=W(z(t),g(t)),\quad t\in I.
 $$
\endproclaim

\head
\dfs{230107-1648} Proof of the main result
\endhead
In order to establish our main result we need uniform estimates for the attractors ${\Cal A}_\eps$
in $H^1_0(\Omega)\times L^2(\Omega)$.
\proclaim{Lemma~\dft{140506-1144}}
Let $f$ be as in Assumption~\rf{140506-1137}. Then there is a constant $C\in \ro0,\infty..$ such that for all $u$, $v\in\R$ and for a.e. $x\in\Omega$, $$|\partial_uf(x,u)|\le C(1+|u|^2),$$ $$|\partial_uf(x,v)-\partial_uf(x,u)|\le C(1+|u|+|v-u|)|v-u|$$ and $$|f(x,v)-f(x,u)-\partial_uf(x,u)(v-u)|\le
C(1+|u|+|v-u|)|v-u|^2.$$

\endproclaim
\proof For all $u$, $v\in\R$ and a.e.  $x\in\Omega$ we have
$$\partial_u f(x,v)-\partial_u f(x,u)=\int_0^1 \partial_{uu}f(x,u+s(v-u))(v-u)\,\D s$$
and
$$f(x,v)-f(x,u)-\partial_uf(x,u)(v-u)=(v-u)^2\int_0^1\theta[\int_0^1 \partial_{uu}f(x,u+r\theta(v-u))\,\D r]\,\D \theta$$
This easily implies the assertions of the lemma.
\eproof

\proclaim{Proposition~\dft{230107-1717}} Let $f$ and $F$ be as in Assumption
\rf{140506-1137}.  Then,
for every measurable function $v\co \Omega\to\R$, both $\hat f(v)$
and $\hat F(v)$ are measurable and for all measurable functions
$u$, $h\co\Omega\to\R$ $$|\hat f(u)|_{L^2}\le |\hat f(0)|_{L^2}+
C(| u|_{L^2}+ |u|^{3}_{L^6}),\leqno\dff{230107-1852} $$ $$|\hat
f(u+h)-\hat f(u)|_{L^2}\le  C| h|_{L^2}+ C(|u|^{2}_
{L^{6}}+|h|^{2}_{L^{6}}) |h|_{L^6}, \leqno\dff{230107-1853}$$
 $$|\hat F(u)|_{L^1}\le
C( |u|^2_{L^2}/2+ |u|^{4}_{L^{4}}/4)+|u|_{L^2}|\hat
f(0)|_{L^2},\leqno\dff{230107-1854}$$ $$
 |\hat F(u+h)-\hat F(u)|_{L^1}\le (|\hat f(0)|_{L^2}+ C(|u|_{L^2}+|h|_{L^2})+
 4C(|u|_{L^{6}}^{3}+ |h|_{L^{6}}^{3}))
 |h|_{L^2},\leqno\dff{230107-1857}$$
and $$|\hat F(u+h)-\hat F(u)-\hat f(u)h|_{L^1}\le\bigl(
C|h|_{L^2}+ C(|u|^{2}_{L^{6}}+ |h|^{2}_{L^{6}})|h|_{L^{6}}\bigr)
|h|_{L^2}. \leqno\dff{230107-1801}$$ Finally, for every $r\in \ro
3,\infty..$ there is a constant $C(r)\in\ro0,\infty..$ such that
for all $u$, $h\in H^1_0(\Omega)$ $$|\hat f(u+h)-\hat
f(u)|_{H^{-1}}\le C(r) |h|_{L^{2}} + C(r)(|u|^{2}_
{L^{6}}+|h|^{2}_{L^{6}}) |h|_{L^{2}}. \leqno\dff{230107-1810}$$
\endproclaim
\proof Lemma~\rf{140506-1144} implies that $f$ satisfies the
hypotheses of \cite{\rfa{PR}, Proposition~3.11}, to which the
reader is referred for details.\eproof

For $s\in \ci2,6..$ we denote by $C_s\in \ro0,\infty..$ an imbedding constant of the inclusion induced map from $H_1$ to $L^s(\Omega)$.
\proclaim{Proposition~\dft{140506-1201}}
Let $f$ be as in  Assumption~\rf{140506-1137}, $I\subset \R$ be an interval, $u$ be a continuous map from $I$ to $H_1$ such that $u$ is continuously differentiable into $H_{0}$ with $v:=\partial(u;H_{0})$.
Then the composite map $\widehat f\circ u\co I\to H_0$ is defined, $\hat f\circ u$ is continuously differentiable into $H_{-1}$ and
$g:=\partial(\widehat f\circ u;H_{-1})=(\widehat{\partial_u f}\circ u)\cdot v$. Moreover, for every $t\in I$,
$$|g(t)|_{H_{-1}}\le C(C_2+C_6|u(t)|_{L^3}^2)|v(t)|_{L^2}\le C(C_2+C_6C_3|u(t)|_{H_1}^2)|v(t)|_{L^2}.\leqno\dff{170506-2121}$$
\endproclaim
\proof It follows from  Proposition \rf{230107-1717} that for
every $w\in H_1$, $\widehat f(w)\in H_0$. Thus $\widehat f\circ u$
is defined as a function from $I$ to $H_0$. Moreover, for every
$t\in I$ and $\zeta\in H_1$, the function $\widehat{\partial_u
f}(u(t))\cdot v(t)\cdot \zeta\co \Omega\to \R$ is measurable and
so by Lemma~\rf{140506-1144} and H\"older's inequality
$$|\widehat{\partial_u f}(u(t))\cdot v(t)\cdot \zeta|_{L^1}\le
C|v(t)|_{L^2}|\zeta|_{L^2}+C|\,|u(t)|^2\,|_{L^3}|v(t)|_{L^2}|\zeta|_{L^6}.$$
It follows that for every $t\in\R$, $g(t)=\widehat{\partial_u
f}(u(t))\cdot v(t)\in H_{-1}$ and~\rf{170506-2121} is satisfied.

Moreover, for $s$, $t\in I$,
$$\aligned&|\widehat{\partial_u f}(u(t))\cdot v(t)-
\widehat{\partial_u f}(u(s))\cdot v(s) |_{H_{-1}}\\&=\sup_{\zeta\in H_{1},\,|\zeta|_{H_1}\le 1}|\widehat{\partial_u f}(u(t))\cdot v(t)\cdot \zeta-
\widehat{\partial_u f}(u(s))\cdot v(s)\cdot \zeta|_{L^1}\\&\le \sup_{\zeta\in H_{1},\,|\zeta|_{H_1}\le 1}T_1(t)(\zeta)+\sup_{\zeta\in H_{1},\,|\zeta|_{H_1}\le 1}T_2(t)(\zeta),\endaligned$$
where
$$T_1(t)(\zeta)=|(\widehat{\partial_u f}(u(t))-
\widehat{\partial_u f}(u(s)))\cdot v(t)\cdot \zeta|_{L^1}$$
and
$$T_2(t)(\zeta)=|\widehat{\partial_u f}(u(s))\cdot( v(t)\cdot \zeta-
 v(s)\cdot \zeta)|_{L^1}$$
By Lemma~\rf{140506-1144}  we obtain, for all $\zeta\in H_1$ with $|\zeta|_1\le 1$,
$$\aligned T_1(t)(\zeta)&\le C|(1+|u(s)|+|u(t)-u(s)|)\cdot|u(t)-u(s)|\cdot\zeta|_{L^2}|v(t)|_{L^2}\\
&\le C|u(t)-u(s)|_{L^3}|v(t)|_{L^2}|\zeta|_{L^6}\\&
+C|u(s)|_{L^6}|u(t)-u(s)|_{L^6}|v(t)|_{L^2}|\zeta|_{L^6}\\&+
C|u(t)-u(s)|_{L^6}|u(t)-u(s)|_{L^6}|v(t)|_{L^2}|\zeta|_{L^6}\\
&\le CC_3C_6|u(t)-u(s)|_{H_1}|v(t)|_{L^2}
+CC_6^3|u(s)|_{H_1}|u(t)-u(s)|_{H_1}|v(t)|_{L^2}\\&+
CC_6^3|u(t)-u(s)|_{H_1}^2|v(t)|_{L^2}
\endaligned$$
and
$$\aligned T_2(t)(\zeta)&\le C|(1+|u(s)|^2)\cdot \zeta|_{L^2}|v(t)-v(s)|_{L^2}\\&\le C(|\zeta|_{L^2}+|\,|u(s)|^2\,|_{L^3}|\zeta|_{L^6})|v(t)-v(s)|_{L^2}
\\&\le C(C_2+C_6^3|u(s)|_{H_1}^2)|v(t)-v(s)|_{L^2}.\endaligned$$
Since $u$ is continuous into $H_1$ and $v$ is continuous into $H_0=L^2(\Omega)$ it follows that
$$\sup_{\zeta\in H_{1},\,|\zeta|_{H_1}\le 1}T_1(t)(\zeta)+\sup_{\zeta\in H_{1},\,|\zeta|_{H_1}\le 1}T_2(t)(\zeta)\to 0\text{ as $t\to s$}$$ so the map
 $(\widehat{\partial_u f}\circ u)\cdot v$ is continuous into $H_{-1}$.

Now, for $s$, $t\in I$, $t\not=s$,
$$\aligned&(t-s)^{-1}|(\widehat{f}\circ u)(t)-(\widehat f\circ u)(s)-\widehat{\partial_u f}(u(s))\cdot v(s)|_{H_{-1}}\\&=\sup_{\zeta\in H_{1},\,|\zeta|_{H_1}\le 1}
(t-s)^{-1}|(\widehat{f}\circ u)(t)\cdot\zeta-(\widehat f\circ u)(s)\cdot\zeta-\widehat{\partial_u f}(u(s))\cdot v(s)\cdot\zeta|_{L^1}\\&\le
(t-s)^{-1}\sup_{\zeta\in H_{1},\,|\zeta|_{H_1}\le 1}T_3(t)(\zeta)+(t-s)^{-1}\sup_{\zeta\in H_{1},\,|\zeta|_{H_1}\le 1}T_4(t)(\zeta)\endaligned$$
where
$$T_3(t)(\zeta)=|g_{t,\zeta}|_{L^1}$$
with $g_{t,\zeta}=(\widehat{f}\circ u)(t)\cdot\zeta-(\widehat f\circ u)(s)\cdot\zeta-\widehat{\partial_u f}(u(s))\cdot (u(t)-u(s))\cdot\zeta$
and
$$T_4(t)(\zeta)=|\widehat{\partial_u f}(u(s))\cdot (u(t)-u(s)-v(s))\cdot\zeta|_{L^1}.$$
Now, by Lemma~\rf{140506-1144}, for all $\zeta\in H_1$ with $|\zeta|_{H_1}\le 1$ and for a.e. $x\in\Omega$
$$|g_{t,\zeta}(x)|\le C(1+|u(s)(x)|+|u(t)(x)-u(s)(x)|)|u(t)(x)-u(s)(x)|^2|\zeta(x)|$$
so
$$\aligned T_3(t)(\zeta)&\le C(|u(t)-u(s)|_{L^3}|u(t)-u(s)|_{L^2}|\zeta|_{L^6})
\\&+
C(|u(s)|_{L^6}|u(t)-u(s)|_{L^6}|u(t)-u(s)|_{L^2}|\zeta|_{L^6})\\&+
C(|u(t)-u(s)|_{L^6}|u(t)-u(s)|_{L^6}|u(t)-u(s)|_{L^2}|\zeta|_{L^6})
\\
&\le CC_6(C_3|u(t)-u(s)|_{H_1}|u(t)-u(s)|_{L^2})
\\&+
CC_6(C_6^2|u(s)|_{H_1}|u(t)-u(s)|_{H_1}|u(t)-u(s)|_{L^2})\\&+
CC_6(C_6^2|u(t)-u(s)|_{H_1}^2|u(t)-u(s)|_{L^2}).\endaligned\leqno\dff{140506-1636}$$
Since $u$ is continuous into $H_1$ and locally Lipschitzian into $H_0=L^2(\Omega)$ it follows from~\rf{140506-1636} that
$$(t-s)^{-1}\sup_{\zeta\in H_{1},\,|\zeta|_{H_1}\le 1}T_3(t)(\zeta)\to 0\text{ as $t\to s$.}$$
We also have
$$\aligned&T_4(t)(\zeta)\le C|(1+|u(s)|^2)\cdot \zeta|_{L^2}|u(t)-u(s)-v(s)|_{L^2}\\&\le
(C|\zeta|_{L^2}+C|\,|u(s)|^2\,|_{L^3}
|\zeta|_{L^6})|u(t)-u(s)-v(s)|_{L^2}\\&\le
C(C_2+CC_6^3|u(s)|_{H_1}^2)
|u(t)-u(s)-v(s)|_{L^2}\endaligned$$
Since $(t-s)^{-1}|u(t)-u(s)-v(s)|_{L^2}\to 0$ as $t\to s$ it follows that $$(t-s)^{-1}\sup_{\zeta\in H_{1},\,|\zeta|_{H_1}\le 1}T_4(t)(\zeta)\to 0\text{ as $t\to s$.}$$
It follows that $\widehat f\circ u$, as a map into $H_{-1}$, is differentiable at $s$ and $\partial_u(\widehat f\circ u;H_{-1})(s)=(\widehat{\partial_uf}\circ u)(s)\cdot v(s)$. The proposition is proved.\eproof
\proclaim{Proposition~\dft{280506-1948}} Let $\eps\in\oi0,\infty..$
be arbitrary.
Define the function $\tilde V=\tilde V_{\eps}\co H_1\times H_{0}\to \R$ by
$$\tilde V(u,v)=(1/2)\langle u,u\rangle_{H_1}+(1/2)\eps\langle v,v\rangle-\int_\Omega F(x,u(x))\,\D x,\quad (u,v)\in H_1\times H_0.
$$
Let $z\co \R\to H_1\times H_0$, $z(t)=(z_1(t),z_2(t))$, $t\in\R$, be a solution of $\pi_\eps$. Then $\tilde V\circ z$ is differentiable
and
$$(\tilde V\circ z)'(t)=-|z_2(t)|^2_{L^2},\quad t\in\R.$$
\endproclaim
\proof This is an application of Theorem \rf{230107-1639} (for the details see \cite{\rfa{PR}, Proposition~4.1}).
\eproof

\proclaim{Proposition~\dft{150506-1201}} Let $\eps\in\oi0,\infty..$
be arbitrary.
Define the function $V=V_{\eps}\co H_0\times H_{-1}\to \R$ by
$$V(v,w)=(1/2)\langle v,v\rangle+(1/2)\eps\langle w,w\rangle_{H_{-1}},\quad (v,w)\in H_0\times H_{-1}.
$$
Let $z\co \R\to H_1\times H_0$, $z(t)=(z_1(t),z_2(t))$, $t\in\R$, be a solution of $\pi_\eps$. Then $z=(z_1,z_2)$ is differentiable as a map into $H_0\times H_{-1}$ with $z_2=\partial(z_1;H_0)$. Let $u=z_1$, $v=z_2$, $w=\partial(v;H_{-1})$ and $g=(\widehat{ \partial_uf}\circ u)\cdot v$. Then the function $\alpha\co\R\to\R$, $t\mapsto V(v(t),w(t))$ is differentiable and for every $t\in\R$
$$\alpha'(t)=-\langle w(t),w(t)\rangle_{H_{-1}}+
\langle g(t),w(t)\rangle_{H_{-1}}.$$

\endproclaim
\proof For $\eps\in\oi0,\infty..$ and $\kappa\in\R$ let ${\bold B}_{\eps,\kappa}\co H_\kappa\times H_{\kappa-1}\to H_{\kappa-1}\times H_{\kappa-2}$
be defined by
$${\bold B}_{\eps,\kappa}(z)=(-z_2,(1/\eps)(z_2+{\bold A}_\kappa z_1)),\quad z=(z_1,z_2)\in H_\kappa\times H_{\kappa-1}.$$
It follows that $-{\bold B}_{\eps,\kappa}$ is $m$-dissipative on $H_{\kappa-1}\times H_{\kappa-2}$
(cf~\cite{\rfa{PR}, proof of Proposition 3.6}).
 Moreover, if $z\co \R\to H_1\times H_0$ is a solution of $\pi_\eps$, then
$$\aligned z(t)&=e^{-{\bold B}_{\eps,2}(t-t_0)}z(t_0)+\int_{t_0}^t(e^{-{\bold B}_{\eps,2}(t-s)}(0,(1/\eps)\widehat f(z_1(s)));H_1\times H_0)\,\D s\\&=
e^{-{\bold B}_{\eps,1}(t-t_0)}z(t_0)+\int_{t_0}^t(e^{-{\bold B}_{\eps,1}(t-s)}(0,(1/\eps)\widehat f(z_1(s)));H_0\times H_{-1})\,\D s,\\&\quad t,t_0\in\R, t_0\le t.\endaligned$$
Since $z(t_0)\in D({\bold B}_{\eps,1})$ and $t\mapsto (0,(1/\eps)\widehat f(z_1(s)))$ is continuous into $D({\bold B}_{\eps,1})$ it follows from~\cite{\rfa{Go}, proof of Theorem~II.1.3 (i)} that $z=(u,v)$ is differentiable as a map into $H_0\times H_{-1}$ with $v=\partial(u;H_0)$.
Now, in $H_{-1}$,
$$\aligned w=\partial(v;H_{-1})&=(1/\eps)(v-{\bold A}_1\circ u+\widehat f\circ u)=
(1/\eps)(v-{\bold A}_0\circ u+\widehat f\circ u).\endaligned$$
It follows from Proposition~\rf{140506-1201} that $w$ is differentiable into $H_{-2}$ and
$$\partial(w;H_{-2})=(1/\eps)(w-{\bold A}_0\circ v+g).$$
Again~\cite{\rfa{Go}, proof of Theorem~II.1.3 (i)} implies that
$$\aligned (v,w)(t)&=e^{-{\bold B}_{\eps,-1}(t-t_0)}(v,w)(t_0)+\int_{t_0}^t(e^{-{\bold B}_{\eps,-1}(t-s)}(0,(1/\eps)g(s));H_{-2}\times H_{-3})\,\D s\\&=
e^{-{\bold B}_{\eps,1}(t-t_0)}(v,w)(t_0)+\int_{t_0}^t(e^{-{\bold B}_{\eps,1}(t-s)}(0,(1/\eps)g(s));H_0\times H_{-1})\,\D s,\\&\quad t,t_0\in\R, t_0\le t.\endaligned\leqno\dff{310506-1808}$$
Now note that the function $V=V_\eps$ is Fr\'echet differentiable and $$DV(v,w)(\widetilde v,\widetilde w)=\langle v,\widetilde v\rangle_{H_0}+\eps\langle w,\widetilde w\rangle_{H_{-1}}.$$
Thus for $(u,v)\in D(-{\bold B}_{\eps,1})=H_1\times H_0$ and $(\widetilde v,\widetilde w)\in H_0\times H_{-1}$
$$\aligned &DV(v,w)(-{\bold B}_{\eps,1}(v,w)+(\widetilde v,\widetilde w))=\langle v,w+\widetilde v\rangle_{H_0}\\&+\eps\langle w,-(1/\eps)w-(1/\eps){\bold A}_1 v+\widetilde w\rangle_{H_{-1}}=\langle v,\widetilde v\rangle_{H_0}-\langle w,w\rangle_{H_{-1}}+\eps \langle w,\widetilde w\rangle_{H_{-1}}.\endaligned$$
Here we have used the fact that
$$\langle w,{\bold A}_1 v\rangle_{H_{-1}}=\langle {\bold A}_1^{-1}w,{\bold A}_1^{-1}{\bold A}_1 v\rangle_{H_{1}}
=\langle {\bold A}_1^{-1}w, v\rangle_{H_{1}}=\langle w,v\rangle_{H_0}$$
as ${\bold A}_1^{-1}w={\bold A}_2^{-1}w\in H_2$.
Defining $W\co (H_0\times H_{-1})\times (H_0\times H_{-1})\to \R$
by
$$W((v,w),(\widetilde v,\widetilde w))=\langle v,\widetilde v\rangle_{H_0}-\langle w,w\rangle_{H_{-1}}+\eps \langle w,\widetilde w\rangle_{H_{-1}}$$
we see that $W$ is continuous. Now it follows from~\rf{310506-1808} and Theorem \rf{230107-1639} that
$\alpha=V_\eps\circ (v,w)$ is differentiable and
$$\alpha'(t)=-\langle w(t),w(t)\rangle_{H_{-1}}+\langle w(t),g(t)\rangle_{H_{-1}},\quad t\in\R.$$
The proof is complete.
\eproof
\proclaim{Proposition~\dft{150506-1330}}
Let $\eps_0\in\oi 0,\infty..$ be arbitrary. Then for every $r\in\ro0,\infty..$ there is a constant $C(r,\eps_0)\in\ro0,\infty..$ such that whenever $\eps\in\lo0,\eps_0..$ and  $z=(u,v)\co \R\to H_1\times H_0$ is a solution of $\pi_\eps$ with $\sup_{t\in\R}( |u(t)|_{H_1}^2+\eps|v(t)|_{H_0}^2)\le r$ and $w:=\partial(v;H_{-1})$, then $$\sup_{t\in\R}(|v(t)|^2_{H_0}+\eps|w(t)|^2_{H_{-1}})\le C(r,\eps_0).$$

\endproclaim
\proof By $C_i(r)\in\ro0,\infty..$, resp. $C_i(r,\eps_0)\in\ro0,\infty..$ we denote various constants depending only on $r$, resp. on $r$ and $\eps_0$, but independent of $\eps\in\lo0,\eps_0..$ and the choice of a solution $z$ of $\pi_\eps$ with $\sup_{t\in\R}( |u(t)|_{H_1}^2+\eps|v(t)|_{H_0}^2)\le r$.

Let $\eps\in\lo0,\eps_0..$ be arbitrary, $\alpha(t)=V_\eps(v(t),w(t))$, $t\in\R$, and $g=(\widehat{\partial_uf}\circ u)\cdot v$. Using~\rf{170506-2121} we see that
$$|g(t)|_{H_{-1}}\le C(1+C_6C_3^2r^2)|v(t)|_{H_0},\quad t\in\R.\leqno\dff{170506-2127}$$
Proposition~\rf{150506-1201} implies that
$$\aligned \alpha'(t)&\le-|w(t)|_{H_{-1}}^2
+(1/2)|g(t)|_{H_{-1}}^2+(1/2)|w(t)|_{H_{-1}}^2
\\&\le
-(1/2)|w(t)|_{H_{-1}}^2+(1/2)
C^2(1+C_6C_3^2r^2)^2|v(t)|_{H_0}^2, \quad t\in \R.\endaligned\leqno\dff{170506-2146}$$
Thus we obtain, for every $k\in\oi0,\infty..$,
$$\aligned&\alpha'(t)+k\alpha(t)\le
(-(1/2)+(k\eps/2))|w(t)|_{H_{-1}}^2\\&+((1/2)
C^2(1+C_6C_3^2r^2)^2+(k/2))|v(t)|_{H_0}^2,\quad t\in\R.\endaligned$$
Choose  $k=k(\eps_0)\in \oi0,\infty..$ such that $(-(1/2)+(k\eps_0/2))< 0$.
Hence we obtain
$$\aligned&\alpha'(t)+k\alpha(t)\le C_1(r,\eps_0)|v(t)|_{H_0}^2
\quad t\in\R.\endaligned$$
Using Propositions~\rf{280506-1948} and~\rf{230107-1717}
we see that
$$\int_{t_0}^t|v(s)|_{H_0}^2\le C_2(r,\eps_0),\quad t,t_0\in\R,\,t_0\le t.$$
It follows that
$$\aligned\alpha(t)&=e^{-k(t-t_0)}\alpha(t_0)+C_1(r,\eps_0)
\int_{t_0}^te^{-k(t-s)}|v(s)|_{H_0}^2\,\D s\\&\le e^{-k(t-t_0)}\alpha(t_0)+C_3(r,\eps_0),\quad t,t_0\in \R,\,t_0\le t.\endaligned\leqno\dff{210506-1323}$$
Using the definition of $\alpha$ we thus obtain from~\rf{210506-1323}
$$\aligned(1/2)|v(t)|^2_{H_0}+(1/2)\eps|w(t)|_{H_{-1}}^2&\le e^{-k(t-t_0)}((1/2)|v(t_0)|^2_{H_0}+(1/2)\eps|w(t_0)|_{H_{-1}}^2)\\&+C_3(r,\eps_0),\quad t,t_0\in\R,\,t_0\le t.\endaligned\leqno\dff{210506-1333}$$
Since for $t\in\R$, $\eps w(t)=-v(t)-{\bold A}_1 u(t)+\hat f(u(t))$ in $H_{-1}$, it follows that
$$\aligned\eps|w(t)|_{H_{-1}}&\le |v(t)|_{H_{-1}}+|u(t)|_{H_{1}}+|\hat f(u(t))|_{H_{-1}}\\&\le
|v(t)|_{H_{-1}}+ C_5(r)\le C_6(r)\eps^{-1/2}+ C_5(r),\quad t\in\R.\endaligned$$
Thus
$$\eps|w(t_0)|_{H_{-1}}^2\le (1/\eps)(C_6(r)\eps^{-1/2}+ C_5(r))^2.\leqno\dff{210506-1340}$$
Furthermore,
$$|v(t_0)|_{H_0}^2\le r/\eps.\leqno\dff{280506-2044}$$
Inserting~\rf{210506-1340} and~\rf{280506-2044} into~\rf{210506-1333} and letting $t_0\to-\infty$ we thus see that
$$|v(t)|^2_{H_0}+\eps|w(t)|_{H_{-1}}^2\le2 C_3(r,\eps_0),\quad t\in\R.$$
This completes the proof.\eproof

 Fix a $C^\infty$-function $\overline\vartheta\co
\R\to \ci 0,1..$ with $\overline\vartheta (s)=0$ for $s\in
\lo-\infty,1..$ and $\overline\vartheta (s)=1$ for $s\in
\ro2,\infty..$. Let $$\vartheta:=\overline\vartheta^2.$$
For $k\in\N$ let the functions
$\overline\vartheta_k\co\R^N\to\R$ and $\vartheta_k\co
\R^N\to\R$ be defined by
$$\overline\vartheta_k(x)=\overline\vartheta(|x|^2/k^2)\text{
and }\vartheta_k(x)=\vartheta(|x|^2/k^2),\quad x\in \R^N.$$
The following theorem (actually a rephrasing of Theorem 4.4 in \cite{\rfa{PR}}) provides the ``tail-estimates''  mentioned in the Introduction:
\proclaim{Theorem~\dft{240107-1223}} Let
Assumptions~\rf{290506-1840} and \rf{140506-1137} be satisfied. Let $\eps_0>0$ be fixed. Choose $\delta$ and
$\nu\in\oi0,\infty..$ with $$\text{
$\nu\le\min(1,\overline\mu/2)$,
$\lambda_1-\delta>0$ and $1-2\delta\eps_0\ge
0$.}$$
  Under these hypotheses, there is a constant $ c' \in\ro0,\infty..$ and for
every $R\in \ro0,\infty..$  there are constants $ M'=
M'(R)$,  $c_k=c_k(R)\in
\ro0,\infty..$, $k\in\N$ with $c_k\to 0$ for $k\to \infty$,
such that for every
  $\tau_0\in
\ro0,\infty..$, every $\eps$, $0<\eps\leq\eps_0$  and every  solution $z(\cdot)$  of $\pi_\eps$
on $I=\ci0,\tau_0..$ with $|z(0)|_Z\le R$
$$|z_1(t)|_{H_1}^2+\eps|z_2(t)^2|_{H_0}\le  c' + M'e^{-2\delta \nu
t},\quad t\in I.$$ If
$|z(t)|_Z\le R$ for $t\in I$, then
$$|\vartheta_kz_1(t)|_{H_1}^2+\eps|\vartheta_kz_2(t)^2|_{H_0}\le c_k + M'e^{-2\delta \nu
t},\quad k\in \N,\, t\in
I.$$ \endproclaim

Now we can prove the following fundamental result:

\proclaim{Theorem~\dft{220506-1206}}
Let $(\eps_n)_n$ be a sequence of positive numbers converging to $0$. For each $n\in\N$ let $z_n=(u_n,v_n)\co \R\to H_1\times H_0$ be a solution of $\pi_{\eps_n}$
such that
$$\sup_{n\in\N}\sup_{t\in\R}( |u_n(t)|_{H_1}^2+\eps_n|v_n(t)|_{H_0}^2)\le r<\infty.$$
Then, for every $\alpha\in\lo0,1..$, a subsequence of $(z_n)_n$ converges in $H_{1}\times H_{-\alpha}$, uniformly on compact subsets of $\R$, to a function $z\co\R\to H_1\times H_0$ with $z=(u,v)$, where $u$ is a solution of $\tilde \pi$ and $v=\partial(u;H_0)$.
\endproclaim
\proof We may assume that $\eps_n\in\lo0,\eps_0..$ for some $\eps_0\in\oi0,\infty..$ and all $n\in\N$. Write $u_n=z_{n,1}$ and $v_n=z_{n,2}$, and $n\in\N$.
We claim that for every $t\in \R$, the set $\{\,u_n(t)\mid n\in\N\,\}$ is relatively compact in $H_0$. Let $\vartheta_k$, $k\in\N$, be as above. Then, choosing $k\in\N$ large enough and using Theorem \rf{240107-1223} we can make $\sup_{n\in \N}|\vartheta _ku_n(t)|_{H_1}$ as small as we wish. Therefore, by a Kuratowski measure of noncompactness argument, we only have to prove that for every $k\in\N$, the set $S_k=\{\,(1-\vartheta_k)u_n(t)\mid n\in\N\,\}$ is relatively compact in $H_0$. Let $U$ be the ball in $\R^N$ with radius $2k$ centered at zero. Then $(1-\vartheta_k)|U\in C^1_0(U)$, so $(1-\vartheta_k)\tilde u_n(t)|U\in H^1_0(U)$ for $n\in\N$. Since $H^1_0(U)$ imbeds compactly in $L^2(U)$ and $(1-\vartheta_k)\tilde u_n(t)|(\R^N\setminus U)\equiv0$, it follows that, indeed, $S_k$ is relatively compact in $H_0$. This proves our claim.

Since, by Proposition~\rf{150506-1330}, for each $n\in\N$, $u_n$ is differentiable into $H_0$  and $v_n=\partial(u_n;H_0)$ is bounded in $H_0$ uniformly $t\in\R$ and $n\in \N$, we may assume, using the above claim and Arzel\`a-Ascoli theorem, and taking subsequences if necessary, that $(u_n)_n$ converges in $H_0$, uniformly on compact subsets of $\R$, to a continuous function $u\co\R\to H_0$. Moreover, since, for each $t\in \R$, $(u_n(t))_n$ has a subsequence that is weakly convergent in $H_1$, we see that $u$ takes its values in $H_1$. Let $w_n=\partial(v;H_{-1})$, $n\in\N$.

For every $n\in\N$ and every $t\in\R$,
$$\eps_n w_n(t)=-v_n(t)-{\bold A}_0u_n(t)+\widehat f(u_n(t))\leqno\dff{300506-1327}$$
in $H_{-1}$.
Now, uniformly for $t$ lying in compact subsets of $\R$, $\widehat f(u_n(t))\to \widehat f(u(t))$ in $H_{-1}$ (by Proposition \rf{230107-1717}), ${\bold A}_0u_n(t)\to {\bold A}_0u(t)$ in $H_{-2}$ and $\eps_n w_n(t)\to 0$ in $H_{-1}$ (by Proposition~\rf{150506-1330}). It follows from~\rf{300506-1327} that, uniformly for $t$ in compact subsets of $\R$, $v_n(t)\to v(t)$ in $H_{-2}$, where $v\co \R\to H_{-2}$ is a continuous map such that, for every $t\in\R$,
$$v(t)=-{\bold A}_0u(t)+\widehat f(u(t)) $$
in $H_{-2}$. It follows that $u$ is differentiable into $H_{-2}$ and $v=\partial(u;H_{-2})$.
Then $u$ is differentiable into $H_{-3}$ and, for all $t\in\R$,
$$\partial(u;H_{-3})(t)=-{\bold A}_{-1}u(t)+\widehat f(u(t))$$
in $H_{-3}$. Since $\widehat f\circ u$ is continuous into $D({\bold A}_{-1})=H_{-1}$ it follows that
$$\aligned u(t)&=e^{-{\bold A}_{-1}(t-t_0)}u(t_0)+\int_{t_0}^t(e^{-{\bold A}_{-1}(t-s)}\widehat f(u(s));H_{-3})\,\D s\\&=
e^{-{\bold A}_{1}(t-t_0)}u(t_0)+\int_{t_0}^t(e^{-{\bold A}_{1}(t-s)}\widehat f(u(s));H_{-1})\,\D s,\quad t,t_0\in\R,\, t_0\le t.\endaligned\leqno\dff{220506-1854}$$
We claim that $u$ is a solution of $\tilde \pi$. To this end let $t_0\in\R$ be arbitrary. Let $\tilde u\co\ro0,\infty..\to H_1$ be the  solution of $\tilde \pi$ with $\tilde u(0)=u(t_0)$ ($\tilde u$ exists by results in~\cite{\rfa{PR3}}). We must show that $\tilde u(s)=u(s+t_0)$ for all $s\in \ro0,\infty..$. If not, then there is a $s_0\ge 0$ with $\tilde u(s_0)=u(s_0+t_0)$ and $\tilde u(s_n)\not=u(s_n+t_0)$ for all $n\in\N$, where $(s_n)_n$ is a sequence with $s_n\to s_0^+$ as $n\to \infty$. By Corollary~\rf{230506-0733} there is a constant $C\in\ro0,\infty..$ such that
$$|e^{-{\bold A}_1 t}w|_{H_0}\le Ct^{-1/2}|w|_{H_{-1}},\quad w\in H_{-1},\,t\in\oi0,\infty...$$
Moreover, by Proposition \rf{230107-1717}, for every $b\in\oi0,\infty..$ there is an $L(b)\in\oi0,\infty..$ such that for all  $w_i\in H_1$,
$|w_i|_{H_1}\le b$, $i=1$, $2$,
$$|\widehat f(w_2)-\widehat f(w_1)|_{H_{-1}}\le L(b)|w_2-w_1|_{H_0}.$$
There is an $\overline s\in\oi s_0,\infty..$ such that whenever $s\in \ci s_0,\overline s..$ then $|u(s+t_0)|_{H_1}<r+1$ and $|\tilde u(s)|_{H_1}<r+1$. Let $L=L(b)$ where $b=r+1$.
Choosing $\overline s$ smaller, if necessary, we can assume that
$$CL(\overline s-s_0)^{1/2}/2<1.\leqno\dff{230506-1027}$$
It follows that, for each $s\in\ci s_0,\overline s..$,
$$u(s+t_0)-\tilde u(s)=\int_{s_0}^s e^{-{\bold A}_1 (s-r)}[\widehat f(u(r+t_0))-\widehat f(\tilde u(r))]\, \D r$$
so
$$\aligned |u(s+t_0)-\tilde u(s)|_{H_{0}}&\le C\int_{s_0}^s (s-r)^{-1/2}L[| u(r+t_0)-\tilde u(r)|_{H_0}]\, \D r\\&\le CL (\overline s-s_0)^{1/2}/2\sup_{r\in\ci s_0,\overline s..}| u(r+t_0)-\tilde u(r)|_{H_0}.\endaligned$$
In view of~\rf{230506-1027}, we obtain that $u(s+t_0)=\tilde u(s)$ for $s\in \ci s_0,\overline s..$, a contradiction, which proves our claim.

We now claim that $u_n(t)\to u(t)$ in $H_1$, uniformly for $t$ lying in compact subsets of $\R$. If this claim is not true, then there is a strictly increasing sequence $(n_k)_n$ in $\N$ and a sequence $(t_k)_k$ in $\R$ converging to some $t_\infty\in\R$ such that
$$\inf_{k\in\N}|u_{n_k}(t_k)-u(t_\infty)|_{H_1}>0.\leqno\dff{310506-0853}$$
For $\eps\in\oi0,\infty..$ define the function $\Cal F_\eps\co H_1\times H_0\to\R$ by
$$\aligned\Cal F_\eps(z)&=(1/2)\eps\langle \delta z_1+z_2,\delta z_1+z_2\rangle+(1/2)\langle A\nabla z_1,\nabla z_1\rangle\\&+(1/2)\langle(\beta-\delta+\delta^2\eps)z_1,z_1\rangle -\int_{\Om}F(x,z_1(x))\,\D x\endaligned$$
where $\delta\in\oi0,\infty..$ is such that $\lambda_1-\delta>0$ and $1-2\delta\eps_0>0$.
Note that
$$\|u\|^2=\langle A\nabla u,\nabla u\rangle+\langle(\beta-\delta)u,u\rangle,\quad u\in H_1$$
defines a norm on $H_1$ equivalent to the usual norm on $H_1$.
Let $\eps\in\lo0,\eps_0..$ and $\zeta=(\zeta_1,\zeta_2)\co\ro0,\infty..\to Z$  be an arbitrary solution of $\pi_\eps$.
Using Theorem \rf{230107-1639} (cf \cite{\rfa{PR}, Proposition~4.1}) one can see that  the function $\Cal F_\eps\circ \zeta$ is continuously differentiable and for every $t\in\ro0,\infty..$
 $$
  \aligned
  &(\Cal F_\eps\circ \zeta)'(t)+2\delta \Cal F_\eps(\zeta(t))=
  \int_{\Om}(2\delta\eps-1)(\delta
  \zeta_1(t)(x)+\zeta_2(t)(x))^2\,\D x\\&+ \int_{\Om}\delta
  \zeta_1(t)(x)f(x,\zeta_1(t)(x))\,\D x -2\delta\int_{\Om}F(x,\zeta_1(t)(x))\,\D x.
  \endaligned\leqno\dff{310506-0914}
 $$
 Moreover, define $\Cal F_0\co H_1\to \R$ by
 $$\aligned
  \Cal F_0(u)&= (1/2)\langle A\nabla u,\nabla u\rangle+(1/2)\langle(\beta-\delta)u,u\rangle-\int _\Omega F(x,u(x))\,\D x
  ,\quad u\in H_1.\endaligned
 $$
 Every solution $\xi\co \R\to H_1$ of $\widetilde \pi$ is differentiable into $H_1$ so  the function $\Cal F_0\circ\xi$ is differentiable and a simple computation shows that
 for $t\in\R$,
 $$\aligned &(\Cal F_0\circ\xi)'(t) +2\delta (\Cal F_0\circ\xi)(t)=-\langle \delta \xi(t)+\eta(t),\delta \xi(t)+\eta(t)\rangle\\&+\int_\Omega [\delta \xi(t)(x)f(x,\xi(t)(x))-2\delta F(x,\xi(t)(x))]\,\D x
 \endaligned\leqno\dff{300506-1320}$$
 where $\eta(t)=-{\bold A}_1\xi(t)+\widehat f(\xi(t))$, $t\in\R$.

Fix $l\in\N$ and, for $k\in\N$,
 let $\zeta_k(t)=z_{n_k}(t_{k}-l+t)$ and
  $\zeta(t)=(u(t_\infty-l+t),v(t_\infty-l+t)$ for $t\in\ro0,\infty..$. Then~\rf{310506-0914} and~\rf{300506-1320} imply that
 $$
  \aligned
   &\Cal F_{\eps_{n_k}}(z_{n_k}(t_k))=e^{-2\delta l}\Cal F_{\eps_{n_k}}(z_{n_k}(t_k-l))\\&+(2\delta\eps_{n_k}-1)\int_0^l e^{-2\delta (l-s)}\left(\int_{\Om}(\delta
  \zeta_{k,1}(s)(x)+\zeta_{k,2}(s)(x))^2\,\D x\right)\,\D s\\&+ \int_0^l e^{-2\delta (l-s)}\left(\int_{\Om}\delta
  \zeta_{k,1}(s)(x)f(x,\zeta_{k,1}(s)(x))\,\D x -2\delta\int_{\Om}F(x,\zeta_{k,1}(s)(x))\,\D x\right)\,\D s.
  \endaligned
  \leqno\dff{070306-0916}
 $$ and
 $$
  \aligned
   &\Cal F_0(u(t_\infty))=e^{-2\delta l}\Cal F_0(u(t_\infty-l))\\&-\int_0^l e^{-2\delta (l-s)}\left(\int_{\Om}(\delta
  \zeta_{1}(s)(x)+\zeta_{2}(s)(x))^2\,\D x\right)\,\D s\\&+ \int_0^l e^{-2\delta (l-s)}\left(\int_{\Om}\delta
  \zeta_{1}(s)(x)f(x,\zeta_{1}(s)(x))\,\D x -2\delta\int_{\Om}F(x,\zeta_{1}(s)(x))\,\D x\right)\,\D s.
  \endaligned
  \leqno\dff{070306-0926}
 $$
 Since $\zeta_{k,1}(s)\to \zeta_{1}(s)$ in $H_0$, uniformly for $s$ lying in compact subsets of $\R$, we obtain from Proposition \rf{230107-1717} that
$$
  \aligned
   &\int_0^l e^{-2\delta (l-s)}\left(\int_{\Om}\delta
  \zeta_{k,1}(s)(x)f(x,\zeta_{k,1}(s)(x))\,\D x -2\delta\int_{\Om}F(x,\zeta_{k,1}(s)(x))\,\D x\right)\,\D s\\&\to
  \int_0^l e^{-2\delta (l-s)}\left(\int_{\Om}\delta
  \zeta_{1}(s)(x)f(x,\zeta_{1}(s)(x))\,\D x -2\delta\int_{\Om}F(x,\zeta_{1}(s)(x))\,\D x\right)\,\D s
  \endaligned
  \leqno\dff{070306-1517}
 $$
 as $k\to\infty$.
 We claim that
 $$
  \aligned
   &\limsup_{k\to \infty}(2\delta\eps_{n_k}-1)\int_0^l e^{-2\delta (l-s)}\left(\int_{\Om}(\delta
   \zeta_{k,1}(s)(x)+\zeta_{k,2}(s)(x))^2\,\D x\right)\,\D s\\&\le
   -\int_0^l e^{-2\delta (l-s)}\left(\int_{\Om}(\delta
  \zeta_{1}(s)(x)+\zeta_{2}(s)(x))^2\,\D x\right)\,\D s.
  \endaligned
  \leqno\dff{070306-1126}
 $$
 In fact, since $1-2\delta \eps_{n_k}\ge0$ for all $k\in\N$
 we have by Fatou's lemma
 $$
  \aligned
   &\limsup_{k\to \infty}(2\delta\eps_{n_k}-1)\int_0^l e^{-2\delta (l-s)}\left(\int_{\Om}(\delta
   \zeta_{k,1}(s)(x)+\zeta_{k,2}(s)(x))^2\,\D x\right)\,\D s\\&=
   -\liminf_{k\to \infty}(1-2\delta\eps_{n_k})\int_0^l e^{-2\delta (l-s)}\left(\int_{\Om}(\delta
   \zeta_{k,1}(s)(x)+\zeta_{k,2}(s)(x))^2\,\D x\right)\,\D s\\&=
   -\liminf_{k\to \infty}\int_0^l e^{-2\delta (l-s)}\left(\int_{\Om}(\delta
   \zeta_{k,1}(s)(x)+\zeta_{k,2}(s)(x))^2\,\D x\right)\,\D s
   \\&\le
   -\int_0^l e^{-2\delta (l-s)}\liminf_{k\to\infty}\left(\int_{\Om}(\delta
   \zeta_{k,1}(s)(x)+\zeta_{k,2}(s)(x))^2\,\D x\right)\,\D s.
  \endaligned
  \leqno\dff{070306-1524}
  $$
  Let $s\in\ci0,l..$ be arbitrary.

   Since $((\zeta_{k,1}(s), \zeta_{k,2}(s)))_k$ converges to $(\zeta_1(s),\zeta_2(s))$ weakly in $H_1\times H_0$  it follows that $((\zeta_{k,1}(s), \delta \zeta_{k,1}(s)+\zeta_{k,2}(s)))_k$ converges to $(\zeta_1(s),\delta\zeta_1(s)+\zeta_2(s))$ weakly in $H_1\times H_0$. It follows that for every $v\in L^2(\Omega)$
  $$\langle v,\delta \zeta_{k,1}(s)+\zeta_{k,2}(s)\rangle\to
  \langle v,\delta \zeta_{1}(s)+\zeta_{2}(s)\rangle\text{ as $k\to\infty$.}$$
  Taking $v=(\delta\zeta_1(s)+\delta\zeta_2(s))$
  we thus obtain
  $$\aligned
  &|(\delta\zeta_1(s)+\delta\zeta_2(s))|_{L^2}^2
  =\langle(\delta\zeta_1(s)+\delta\zeta_2(s)),
  (\delta\zeta_1(s)+\delta\zeta_2(s))\rangle\\&=
  \lim_{k\to\infty}\langle (\delta\zeta_1(s)+\delta\zeta_2(s)),
  (\delta\zeta_{k,1}(s)+\delta\zeta_{k,2}(s))\rangle\\&
  \le| (\delta\zeta_1(s)+\delta\zeta_2(s))|_{L^2}\liminf_{k\to\infty}| (\delta\zeta_{k,1}(s)+\delta\zeta_{k,2}(s))|_{L^2}
  \endaligned
   $$
  and so
  $$\aligned
   &\int_{\Om}(\delta
   \zeta_{1}(s)(x)+\zeta_{2}(s)(x))^2\,\D x\le
   \liminf_{k\to\infty}\int_{\Om}(\delta
   \zeta_{k,1}(s)(x)+\zeta_{k,2}(s)(x))^2\,\D x.
      \endaligned
      \leqno\dff{070306-1551}
  $$
  Inequalities~\rf{070306-1551} and~\rf{070306-1524} prove~\rf{070306-1126}.
  Since, by Proposition \rf{230107-1717},
  $$\int_{\Omega}F(x,u_{n_k}(t_k)(x))\,\D x\to \int_\Omega F(x,u(t_\infty)(x))\,\D x$$ we obtain,
 using Proposition~\rf{150506-1330}, that
 $$\limsup_{k\to\infty}\Cal F_{\eps_{n_k}}(z_{n_k}(t_k))=(1/2)\limsup_{k\to\infty}\|u(t_k)\|^2-
 \int_{\Omega}F(x,u(t_\infty)(x))\,\D x$$
 Moreover, there is a constant $C'\in\oi0,\infty..$ such that
 $$\sup_{k\in\N}\sup_{t\in\R}|\Cal F_{\eps_{n_k}}(z_{n_k}(t))|+\sup_{t\in\R}|\Cal F_{0}(u(t))|\le C'.$$
 Thus
 $$\aligned
  &(1/2)\limsup_{k\to\infty}\|u(t_k)\|^2-
 \int_{\Omega}F(x,u(t_\infty)(x))\,\D x\le
  e^{-2\delta l}C'\\&
  -\int_0^l e^{-2\delta (l-s)}\left(\int_{\Om}(\delta
  \zeta_{1}(s)(x)+\zeta_{2}(s)(x))^2\,\D x\right)\,\D s\\&+ \int_0^l e^{-2\delta (l-s)}\left(\int_{\Om}\delta
  \zeta_{1}(s)(x)f(x,\zeta_{1}(s)(x))\,\D x -2\delta\int_{\Om}F(x,\zeta_{1}(s)(x))\,\D x\right)\,\D s\\&=
  e^{-2\delta l}C' +(1/2)\|u(t_\infty)\|^2-\int_{\Omega}F(x,u(t_\infty)(x))\,\D x\\&-e^{-2\delta l}\Cal F_0(u(t_\infty-l))\le 2e^{-2\delta l}C'+(1/2)\|u(t_\infty)\|^2-\int_{\Omega}F(x,u(t_\infty)(x))\,\D x.
 \endaligned
 $$
 Thus for every $l\in\N$
 $$
\limsup_{k\to\infty}\|u(t_k)\|^2\le 4e^{-2\delta l}C'+\|u(t_\infty)\|^2
 $$
 so
 $$
\limsup_{k\to\infty}\|u(t_k)\|\le \|u(t_\infty)\|.
 $$
 Since $(u_{n_k}(t_{n_k}))_k$ converges to $u(t_\infty)$ weakly in $H_1$ we have
 $$
  \liminf_{k\to\infty}\|u_{n_k}(t_{n_k})\|\ge\|u(t_\infty)\|.
 $$
 Altogether we obtain
 $$
  \lim_{k\to\infty}\|u_{n_k}(t_{n_k})\|=\|u(t_\infty)\|.
 $$
 This implies that $(u_{n_k}(t_{n_k}))_k$ converges to $u(t_\infty)$ strongly in $H_1$, a contradiction to~\rf{310506-0853}. Thus, indeed, $u_n(t)\to u(t)$ in $H_1$, uniformly for $t$ lying in compact subsets of $\R$.

Now~\rf{300506-1327} implies that $v_n(t)\to v(t)$ in $H_{-1}$, uniformly for $t$ lying in compact subsets of $\R$. Since $(v_n)_n$ is bounded in $H_0$, interpolation between $H_0$ and $H_{-1}$ (cf. \cite{\rfa{PR3}}) now implies that $v_n(t)\to v(t)$ in $H_{-\alpha}$, uniformly for $t$ lying in compact subsets of $\R$.
 The proof is complete.
\eproof
Now we obtain the main result of this paper.
\proclaim{Theorem~\dft{151106-1022}} For every $\alpha\in\lo0,1..$
the family $(\Cal A_\eps)_{\eps\in\ro0,\infty..}$ is upper semicontinuous at $\eps=0$ with respect to the topology of $H_1\times H_{-\alpha}$, i.e.
$$\lim_{\eps\to 0^+}\sup_{y\in\Cal A_\eps}\inf_{z\in \Cal A_0}|y-z|_{H_1\times H_{-\alpha}}=0.$$
\endproclaim
\proof
Using the first part of Theorem \rf{240107-1223}, choosing $\eps_0\in\oi0,\infty..$ arbitrarily and $\delta\in\oi0,\infty..$ such that $\lambda_1-\delta>0$ and $1-2\delta\eps_0>0$ and noting that the constant $c'$ in that theorem is independent of $\eps\in\lo0,\eps_0..$, it follows that for all $\eps\in \lo0,\eps_0..$ and all $(u,v)\in \Cal A_\eps$,
$$|u|^2_{H_1}+\eps|v|^2_{H_0}\le 2c'.$$
Now an obvious contradiction argument using Theorem~\rf{220506-1206} completes the proof of our main result.\eproof
\remark{Remark} Theorem~\rf{151106-1022} and Corollary~\rf{300506-1029} imply Theorem~\rf{300506-0938}.\endremark

\enddocument